\documentclass[12pt,titlepage,a4paper
,draft
]{article}
\usepackage{amssymb,latexsym,amsfonts,amsmath,amsthm}

\usepackage{verbatim}

{\makeatletter \@addtoreset{equation}{section}}

\newcommand{\diag}{\mathop{\mathrm{diag}}\nolimits}
\newcommand{\Ker}{\mathop{\mathrm{Ker}}\nolimits}
\newcommand{\IIm}{\mathop{\mathrm{Im}}\nolimits}

\newcommand{\C}{{\mathbb{C}}}
\newcommand{\R}{{\mathbb{R}}}

\newcommand{\h}{ \mathop{ \mathrm{h} {} }\nolimits }
\newcommand{\e}{ \mathop{ \mathrm{e} {} }\nolimits }

\newcommand{\beq}{\begin{equation}}
\newcommand{\eeq}{\end{equation}}

\newtheorem{theorem}{Theorem}[section]

\newtheorem{corollary}[theorem]{Corollary}
\newtheorem{prop}[theorem]{Proposition}

\begin{document}

\title{A class of marked invariant subspaces with an application to
algebraic Riccati equations}

\author{Pudji Astuti
\\
Faculty of Mathematics and \\
 Natural Sciences\\
Institut Teknologi Bandung\\
Bandung 40132\\
Indonesia
         \and
Harald~K. Wimmer\\
Mathematisches Institut\\
Universit\"at W\"urzburg\\
D-97074 W\"urzburg\\
Germany}
%


\date{}

\maketitle

\begin{abstract}

\vspace{2cm}

\vspace{.2cm}



\vspace{.5cm}

\noindent
{\bf Abstract:\quad}  Invariant subspaces of a matrix $A$ 
are considered which are obtained by truncation of a
 Jordan basis 
of a generalized eigenspace of $A$. We characterize those 
subspaces which are independent of the choice of the
Jordan basis. An application to Hamilton matrices
and algebraic Riccati equations is given.



\end{abstract}



\section{Invariant subspaces}

Let  $ \lambda $ be an  eigenvalue of a  complex $n \times n$ matrix 
 $A$ and let 
\[ 
E _{  \lambda }(A)  = \Ker ( A -  \lambda I) ^n
\] 
be the corresponding  generalized eigenspace.  
Suppose \,
$ \dim  E _{  \lambda }(A) = k $. If 
\[
( s -   \lambda ) ^{t_1}, \, \dots \, , ( s -   \lambda ) ^{t_k},  \,\, \,
  t_1 \leq \cdots  \leq   t_k ,
\]
are the corresponding elementary divisors 
then  \,$E _{  \lambda }(A) $ 
is a direct sum of  $t_i$-dimensional
cyclic subspaces,
i.e.  
\[
  E _{  \lambda }(A) = K_1 \oplus \cdots  \oplus  K_k
\]
with 
\beq \label{eq.jch}
K_i 
=  
 {\rm{span}} \bigl\{ u_i, \,  ( A -  \lambda I)u_i,
 \, \dots ,   ( A -  \lambda I) ^{ t_i -1 } u_i  \bigr\} ,
\eeq 
and \,$  ( A -  \lambda I) ^{ t_i} u_i  = 0$, \, 
$  i= 1,\dots , k$ .
We call \beq  \label{eq.uu}
U = ( u_1, \dots , u_k )
\eeq 
 a tuple of \emph{generators} of 
$  E _{  \lambda }(A)  $.
From a given   $U$  one can construct 
$A$-invariant subspaces 
in the following way.
Let  $r = (r_1, \dots , r_k )$ be such that
\beq \label{eq.co1} 
 0 \,  \leq  \,  r_i   \, <   \, t_i , \, \,   i = 1, \dots ,k .
\eeq
We set 
\beq  \label{eq.ri}
     W_{r_i} (U)  = {\rm{span}}   \bigl\{ 
 ( A -  \lambda I)^{r_i}u_i, \, ( A -  \lambda I) ^{r_i +1}u_i, \, 
 \dots ,  ( A -  \lambda I) ^{t_i -1} u_i \bigr\}
\eeq 
and 
\beq \label{eq.dis}
 W(r,U)  \,  =  \, 
           W_{r_1}(U) \, \oplus  \,  \cdots   \,  \oplus   \, W_{r_k}(U) .
\eeq
The construction of invariant subspaces
of the from $W(r,U)$
is a standard procedure in linear algebra and systems theory
(see e.g. \cite{Ku}, \cite[p.61]{GLR0}, \cite{Fr}, \cite{RR},
\cite[p.28]{BEGO}).

If $U$ and $\tilde{U}$ are two different tuples of generators of
$  E _{  \lambda }(A)  $ then the restrictions of $A$ to  $W(r,U)$ 
and $ W(r, \tilde{U})$ have the same elementary divisors,
namely  
$ (s-  \lambda) ^{t_i-   r_i}$, $ i = 1, \dots ,k$. 
However,
 in general, the subspaces 
$ W(r,U) $ and $ W(r, \tilde{U})$ will be different. 
  Consider the following  example with
$k= 2$, $t_1= 2$, $t_2 = 3$, and 
\beq \label{eq.mtr}
 A = \diag (N_2, N_3), 
\, 
          N_2 = \begin{pmatrix}   0 & 1 \\  0 & 0 
\end{pmatrix}\!,
 \, 
 N_3 = \begin{pmatrix}   0 & 1   & 0 
\\ 
 0 & 0  & 1 \\
                0 & 0 & 0
\end{pmatrix}\!.
\eeq 
Let   $e_i$ 
be a unit vector of $\C^5$.
Then 
  $U =  \{ e_2 , e_5 \}$ and
$ \tilde{U} =  \{ e_2, e_5 + e_2 \} $ are tuples of generators of
$E_0(A) = \Ker A^5 = \C^5$. 
If we choose  \,$ r = (1,0) $, then 
$  W(r,U) 
  = {\rm{span}} \{ e_1, e_3, e_4, e_5  \} $
  and $  W(r, \tilde{U})
  = 
 {\rm{span}} \{e_1,  e_3 ,  e_4 + e_1, e_5 + e_2  \}$.
Thus  
\beq \label{eq.uuti}
W(r,U) \neq  W(r, \tilde{U}). 
\eeq 
   On the other hand, if we choose \,$ r = (1,2) $,  then   
\beq \label{eq.sm}
W(r,U) =   W(r, \tilde{U}).
\eeq

It is the purpose of our note to determine those 
 tuples $r = ( r_1, \dots , r_k)$
 which have the property that the space 
$ W(r,U) $ given by \eqref{eq.ri} and  \eqref{eq.dis}
is independent of the generator tuple $U$.
The motivation for our study comes from 
Kucera's survey article \cite{Ku}, which
deals with  independence of generator tuples in
the case of  Hamiltonian matrices. 
In  Section \ref{sec.app} 
we make the connection with  \cite[p.60]{Ku} 
applying a corollary of our main theorem
to Hamiltonian matrices and algebraic Riccati
equations. 

In the sequel we  
 assume that $\lambda =  0 $ 
is an eigenvalue of $A$ and we 
focus on $E_0(A) = \Ker A ^n $. 
With  each  nonzero vector  $v \in  E _{0}(A)$  
we associate a {\emph{height}} $\h(v)$  and  an {\emph{exponent}} $\e(v)$
as follows. 
Suppose 
\[
 v \in \IIm A^q, \,   v \notin 
\IIm A ^{q +1},  \,\, 
 v \in \Ker  A  ^p, \, 
v \notin  \Ker  A  ^{p-1}.
\]
Then we  set    
$\h(v) = q$ and $ \e(v) = p$. 
Thus, if $ \lambda = 0 $ in \eqref{eq.jch}
 then the elements of $U$ in \eqref{eq.uu} 
satisfy 
$\e(u_1) = t_1 \leq \cdots \leq \e(u_k) = t_k$\, 
and $\h(u_i) = 0$. 
We define
\[
   \langle v \rangle =  {\rm{span}} \{ A ^{ \nu } v,  \nu \geq 0 \}.
\]
Then $ \langle v \rangle $
is a cyclic subspace generated by
$v$, and  \,$\dim \langle v \rangle =  \e(v) $.


\section{The main result}

\begin{theorem} \label{thm.mn}
Let  
$A \in \C^{n \times n} $  
and let 
\beq \label{eq.elth}
  s^{ t_1 }, \dots ,  s^{ t_k }, \, 
t_1 \leq   \cdots \leq  t_k ,
 \eeq 
be the 
elementary divisors corresponding to the eigenvalue $\lambda = 0$. 
 Let  
\[
   U =   (u_1,  \dots , u_k )
\]
be a tuple  of  generators of 
$E_0(A) = \Ker A^n$
such that 
\, $
 \e(u_i) = t_i$, $i = 1, \dots , k$, 
and 
\[
E_0(A) =  \langle u_1 \rangle \, \oplus \, \cdots \,  
\oplus \, 
    \langle u_k\rangle .
\]
Let  \,$ r = ( r_1, \dots , r_k) $ be a $k$-tuple of integers
with \,$ 0 \leq r_i < t_i$,   $i = 1, \dots , k$.
Define  
\beq  \label{eq.agn}
 W(r,U) 
              = 
  \langle  A ^{r_1} u_1 \rangle \, \oplus \, \cdots \,  
\oplus \, 
    \langle  A ^{r_k}  u_k \rangle
\eeq
and 
\beq \label{eq.wr}
  W(r) = \bigl(\IIm A ^{r_1} \cap \Ker A^{t_1 - r_1} \bigr)
+
\cdots + \bigl(\IIm A ^{r_k} \cap \Ker A^{t_k - r_k} \bigr).
\eeq
Then the following statements are equivalent:
\\
{\rm{(i)}} 
The $k$-tuple $r = ( r_1,   \cdots , r_k) $
satisfies 
\beq \label{eq.rmo}
 r_1  \leq \cdots  \leq r_k,
\eeq
and
 \beq \label{eq.df}
t_1 - r_1 \leq \cdots \leq t _k  - r_k.
\eeq
{\rm{(ii)}}
The space
  $  W(r,U) $
is independent of  \,$U$.\\
Moreover, if \eqref{eq.rmo} and  \eqref{eq.df} hold then  
 \,$W(r,U) = W(r)$.
\end{theorem}

\medskip 
\noindent Proof. 
(i) $\Rightarrow$ (ii). 
We show that 
  \eqref{eq.rmo} and \eqref{eq.df}
imply 
$ W(r,U) = W(r)$.  
Define
\,$
   W_{r_s}(U) = \langle  A ^ {  r_s } u_s \rangle $
such that \eqref{eq.dis} holds.
From 
\[
   W_{r_s} (U) \subseteq 
 \IIm A ^{r_s}\, \cap  \, \Ker  A^{t_s - r_s} 
\] 
 we immediately obtain  $  W(r,U)  \subseteq  W(r) $. 
Now let 
 $ x $ be in  $\IIm A ^{r_s} \cap \Ker A^{t_s - r_s}  $. 
Then 
\, $ x =  A ^{r_s} y$ \, for some \, $ y \in E_0(A)$, 
and 
\beq \label{eq.zro} 
  A^{t_s - r_s} x = A^ {t_s} y  = 0. 
\eeq 
With respect to the basis  
\beq  \label{eq.mclb}
  \mathcal{B} _U 
= 
 \{A^{\nu_i}u_i; \,  0 \leq \nu _i \leq t_i -1 , \,
i = i, \dots , k \} 
\eeq
 we have
\[   y  
= 
  \sum _{i = 1}^ k \,  \sum _{\nu_i  = 0} ^{t_i -1}
\,
  \alpha _{i \nu _i} A^{ \nu_i } u_i .
\] 
Let $\ell$ be the 
largest integer sucht that  \,$t_{ \ell } \leq t_s$.
Then $ A^{ t_s} u_i = 0$ for \, $i = 1, \dots, 
\mbox{$\ell$}$.
Moreover 
 $ A^{  t_s + \nu  _i }  u_i = 0 $ if $  t_s + \nu  _i >  
t_i $.
Therefore 
\[ 
 A^{ t_s} y  = 
  \sum _{i > \ell} 
\,  \sum _{\nu _i  = 0}  ^{ t_i - t_s  -1}  
\,
  \alpha _{i \nu _i } A^{  t_s + \nu  _i } u_i  = 0.
\] 
Since the vectors of  $\mathcal{B}_U$ 
are linearly independent we obtain
\,  $\alpha _{i \nu _i} = 0$ for $i > \ell$
and 
\, $ \nu _i = 0, \dots ,  t_i - t_s -1 $. Hence
\[
  y 
=  
\sum _{i = 1}^ {\ell } \,  \sum _{\nu _i  = 0} 
^{ t_i  -1}
             \alpha _{i \nu _i }
                          A^{  \nu  _i } u_i \,
+ \, 
  \sum _{i > \ell} 
                       \sum _{\nu _i  = t_i - t_s }  ^{t_i -1}
  \alpha _{i \nu _i }  
                         A^{  \nu  _i } u_i
\]
and        
\[
x
  = 
  \sum _{i = 1}^ {\ell} \,  \sum _{\nu _i  = 0} 
^{ t_i  -1}
             \alpha _{i \nu _i }
                          A^{ r_s + \nu  _i } u_i \,
+ \, 
  \sum _{i > \ell} 
                       \sum _{\nu _i  = t_i - t_s }  ^{t_i -1}
  \alpha _{i \nu _i }  
                         A^{r_s +  \nu  _i } u_i.  
\]
  Note that  
$   t_s = \cdots = t_{\ell} $ 
 implies  $ r_s  = \cdots = r_{\ell}.$
Hence, 
if \,$ 1 \leq i \leq  \ell$ then  
$ r_i \leq  r_s $, and therefore 
    \beq \label{eq.dp} 
  A^{ r_s + \nu  _i } u_i \in  W_{r_i} (U).
\eeq
 On the other hand, if 
\,
  $i > \ell $ %
then
      $t_s - r_s \leq t_i - r_i$.
In that case
 $   \nu  _i \in  
 \{  t_i - t_s, \dots ,  t_i  -1 \}  $
implies  
\[
  r_s +  \nu  _i \geq  r_s + (t_i - t_s) \geq r_i. 
\]
Thus,  we again have \eqref{eq.dp}. 
Hence $x \in W(r,U)$ and therefore $ W(r) \subseteq  W(r,U)$.

\medskip 

\noindent 
(ii) $\Rightarrow$ (i). 
We assume that $ W(r,U)$ is    independent of $U$. 
Let us show first that
\beq \label{eq.co2}
 r_i = r_j \,\,\,  {\rm{if}} \,\,\,  t_i = t_j.
 \eeq
Suppose \,$r = (r_1, \dots , r_k) $ is such that
\,$ t_s = t _{s+1}$ and
\,$r_{s} \neq r _{s+1}$,
e.g.
\beq \label{eq.sptt}
 r_{s+1} < r _{s} \,\,\, {\rm{for \,\, some}} \,\,\,
 s  \in \{1, \dots , k-1 \}.
\eeq
 Let
\,$ V = (v_1,   \dots , v_k) $
be such that
 \,$(v_ {s}, v _{s+1 } )   = (u _{s+1 },  u_{s})$,
and \,$v_i = u_i$ if $i \notin \{s, s +1 \}$.
Then
 $A ^{ r_{s+1 } } u_{s+1} \in W(r,U)$ 
 but
 $ A ^{ r_{s+1 } } u_{s+1} =  A ^{ r_{s+1 } }v_{s}  \notin  W(r, V)$.
Therefore the tuples
$U $ and $V$ contain the same elements, but
$ W(r,U) \neq  W(r, V )$.

Now suppose that  \eqref{eq.rmo} is not satisfied.
Then
we have \eqref{eq.sptt},   and
\beq \label{eq.wzv}
A ^ {  r_{s+1} }  u _{s} \notin W(r,U) .
\eeq
Let    $V = (v_1, \dots , v_k )$
be given by
\, $ v_{s+1} = u _{s+1} +  u _{s} $, and $ v_i = u_i $, if
$i \neq s+1$.
Thus $V$ is
  a tuple of generators of $ E_{0}(A) $ with
$\e (v_i) = \e (u _i)$.
Consider
\[
   A ^ {  r_{s+1} }  u _{s+1}
+
A ^ {  r_{s+1} }  u _{s } \, = 
 \, A ^ {  r_{s+1} }  v_{s+1} \,  \in \,  W(r,V).
\]
Then  $   A ^ {  r_{s+1} }  v_{s+1} \notin  W(r,U) $.
Otherwise   $A ^ {  r_{s+1} }  u_{s+1} \in  W(r,U) $
   would imply  \,
$ A ^ {  r_{s+1} }  u _{s} \in  W(r,U) $, which is a
contradiction
to \eqref{eq.wzv}.

\medskip 

Suppose    $r = (r_1, \dots , r_k)$
does not satisfy  \eqref{eq.df}.
Then
 \, $ t_s  - r _s >  t_{s+1}  - r _{s +1} $ \,
for some $s \in \{1, \dots , k -1 \} $.
 Because of \eqref{eq.co2} we
have $t_{s+1} \neq t_s$. Hence
\, $
 r _ {s+1} -  r _s >  t_{s+1} - t_s  > 0 $,
and  \, $  r _s <  r _ {s+1}  $, and
\beq \label{eq.dec}
    r _s + (  t_{s+1} - t_s ) <  r _ {s+1}.
\eeq
Because  \eqref{eq.agn}  it is obvious
that \eqref{eq.dec} implies
\beq \label{eq.spt}
  A ^ {   r _{s}  + (  t_{s+1} - t_{s} ) } u_{s+1} \notin  W (r, U) .
\eeq
Define
\,$
v_s = u _s + A ^{ t_{s+1} -  t_s } u _ {s+1}$.
Then \,$\e(v_s) = \e(u_s) = t_s$.
Therefore
\beq \label{eq.vv}
V = \{ u_1, \dots, u_{s-1}, v_s,  u_{s+1}, \dots, u_k   \}
\eeq
is  another tuple of generators of $E_0(A)$.
Let us show that $ W (r, V) \neq W (r, U) $. Clearly, the vector
$ A ^{r_s} v_s $ belongs to  $ W (r, V) $. Suppose \,
\[
 A ^{r_s} u_s + A ^{ r_s + (  t_{s+1} -  t_s ) } u _ {s+1}
\,
  =
\,
   A ^{r_s} v_s  \,  \in  \,  W (r, U) .
\]
Because of $  A ^{r_s} u_s  \in  W (r, U) $
that would imply
\[
   A ^ {r_s + (t_{s+1} -  t_s)  }  u _ {s+1} \in  W (r, U) ,
\]
which is a contradiction to \eqref{eq.spt}.
\hfill  $\square$ 

\medskip 

\bigskip 
Let us consider again Example \eqref{eq.mtr}.
We have  $(t_1,t_2) = (2,3)$.
In the case of  $r = (1,0)$ condition \eqref{eq.rmo} is violated,
which accounts for
\eqref{eq.uuti}.  
In the case of  $r = (1,2)$
both  \eqref{eq.rmo} and \eqref{eq.df} hold, which ensures  
\eqref{eq.sm}.


In accordance with  a definition in
 \cite[p.\,83]{GLR}  and  \cite{Bru} the space $ W(r,U) $ is a 
  \emph{marked}  $A$-invariant subspace of $E_0(A)$.
That means
 $ W(r,U) $ 
 has a Jordan basis, in our case
\[
 \{A^{r_i+ \mu_i}u_i; \,  0 \leq \mu _i \leq t_i -r_i -1, \,
i = i, \dots , k \}, 
\]
which can be extended to a Jordan basis of  $E_0(A)$,
namely to 
 $\mathcal{B}_U$ in \eqref{eq.mclb}.
Let 
\,$  \mathcal{M}_r $ be the set of marked subspaces $M$ of  $E_0(A)$ 
such  that 
 the elementary divisors of
the restriction $ A _{\mid M} $  are 
$ s^{t_1 - r_1}, \dots ,  s^{t_k - r_k} $.
 We have noted before  that for each tuple of generators  $U$
the corresponding space 
$   W(r,U) $ is in  $ \mathcal{M}_r $.
Suppose \eqref{eq.rmo} and  \eqref{eq.df}
hold. Then all the spaces $   W(r,U) $ coincide with
$W(r)$ and one might ask whether   $W(r)$
is the only  subspaces in $  \mathcal{M}_r $. 
In  the following  we have an example
where $ \mathcal{M}_r \supsetneqq \{  W(r) \}$.
  Let $n = 10$, $k = 2$,
 and $ t = (t_1,t_2 ) = (4,6)$,  and $ r = ( 2, 3) $.
Then $ t - r = ( 2,  3)$. Hence 
the conditions  \eqref{eq.rmo} and  \eqref{eq.df}
are satisfied. 
Let 
$U = ( u_1, u_2) $ be a tuple of generators
such that 
$ \e(  u_1 ) =  4$ and $ \e(  u_2 ) = 6$. 
The subspaces 
\[
  M = W(r,U) = W(r)  = \langle A^2  u_1 \rangle 
\oplus  \langle A^3  u_2   \rangle
\]
and 
\,$
   \tilde{M } = 
 \langle A  u_1 \rangle 
\oplus  \langle A^4  u_2   \rangle
$\, 
are marked,
the elementary divisors  of $A_{\mid M } $ and  
 $A_{\mid  \tilde{M } } $ are $s^2, s^3$.
Hence $  \tilde{M } \in  \mathcal{M}_r $, but  
$  \tilde{M } \neq W(r)$.


\medskip 

Let $[m]$ denote the greatest integer of $m$. 
If we assume  $(t_1, \dots, t_k) $ as in \eqref{eq.elth} and take 
$r = ( [\tfrac{1}{2}t_1  ], \dots ,  [\tfrac{1}{2}t_ k ] )$
 then the conditions  \eqref{eq.rmo} and \eqref{eq.df}
are satisfied and we note the 
following corollary of Theorem \ref{thm.mn}. 
\begin{corollary} \label{eq.spc}
Let  
 $A \in \C^{n \times n}$ and $ 0 \in \sigma(A) $.
Let \,$
  s^{2m_1}, \dots , s^{2m_k}$, 
be the elementary 
 divisors of  $A$ 
 corresponding to $\lambda = 0$.
If 
$ U = ( u_1, \dots , u_k ) $ is  a tuple of generators of \,
$\Ker A^n$
 such that $\e(u_i) = 2m_i$,
$i = 1, \dots , k$,
then  $ \e( A^{m_i}u_i) = m_i$ \, 
for all $i$, and 
\begin{multline*}
\langle  A ^{m_1} u_1 \rangle \oplus 
\cdots   \oplus 
 \langle     A ^{m_k} u_k    \rangle 
= 
\bigl(\IIm A ^{m_1} \cap \Ker A ^{m_1}\bigr)
+
\cdots +
 \bigl( \IIm A ^{m_k} \cap \Ker A ^{m_k}\bigr).
\end{multline*} 
\end{corollary}


\section{An application} \label{sec.app}
In this section 
we apply 
Corollary~\ref{eq.spc} to the algebraic Riccati equation
\beq \label{eq.are}
  Q +   F^* X + XF  -  XDX = 0
\eeq
and its
associated
Hamiltonian matrix
\beq \label{eq.ham}
   H =    \begin{pmatrix}
                              F & - D \\  -Q & - F^*
\end{pmatrix}. 
\eeq 
Here $F, D , Q $ are complex $m \times m$ matrices,
$D$ and $Q$ are hermitian, $D \geq 0$, and the pair 
$(F, D) $ is assumed to be controllable.
Then (see \cite[p.59]{Ku}
all elementary divisors corresponding to eigenvalues
$ i \alpha  \in i \R $ have even degree. 
To fix ideas we assume 
\, $ \sigma (H) = \{ 0\} $.
The subsequent result complements  Lemma 3.2.3 of
\cite[p.60]{Ku}.
\begin{prop}
Let 
 \,$   s ^{2m_1 } ,   \dots  s   ^{ 2m_k }$
be the elementary divisors of $H$. 
Set 
\beq \label{eq.een}
W =  \bigl(\IIm H ^{m_1} \cap \Ker H^{m_1} \bigr)
+
\cdots + \bigl(\IIm H ^{m_k} \cap \Ker H^{m_k} \bigr).
\eeq 
Then $W$ is an $H$-invariant subspace of
$ \C^{2m} $ and 
                $ \dim W = m$.
Let   $ Y , Z \in \C^{m \times m}$  be such that
   the columns of \,$ \left( \begin{smallmatrix}  Y \\ Z 
    \end{smallmatrix}  \right) $ 
are a basis of   $W$.
Then 
$ Y $ is nonsingular and
$X = Z Y ^{-1} $ is the unique hermitian solution  of
\eqref{eq.are}. 
\end{prop} 

Proof.  Set $ t  = (2 m _1, \dots , 2 m_k) $. 
Let $ U = (u_1, \dots , u_k) $ be a tuple of
generators of $E_0(H) = \C ^{2m}$. 
According to  \cite{Ku} we have
\[
W ( \tfrac{1}{2} t , U) = {\rm{span}} \left( \begin{matrix}  I_m \\ X
    \end{matrix}  \right),
\]
where 
$X \in \C^{m \times m}$ is the unique hermitian 
solution of \eqref{eq.are}. 
From Corollary~\ref{eq.spc} we know that
$ W ( \tfrac{1}{2} t, U) $ is independent of the choice of $U$.
Moreover,  
$  W ( \tfrac{1}{2} t, U) = W $ where 
$W $ is given by \eqref{eq.een}.
Hence, if \,$  W = {\rm{span}} \left(  \begin{smallmatrix}
                              Y \\ Z 
 \end{smallmatrix} \right) $
then  $ Y $ is nonsingular, and
\[       
 {\rm{span}}  
  \left(  \begin{matrix}
                              Y \\ Z 
  \end{matrix} \right) 
=  
  {\rm{span}}     \left(  \begin{matrix}
                              I \\ Z  Y ^{-1} 
  \end{matrix} \right) 
\]
implies that  $ X =  Z  Y ^{-1}  $ is the solution of 
\eqref{eq.are}.
\hfill $\square$

\section{Conclusions}
Results of this note can 
be considered in a module theoretic framework. 
In a subsequent paper we shall make the connection
of  Theorem~\ref{thm.mn} 
with marked subspaces in \cite{FPP} and with 
 torsion modules over discrete valuation domains
  in~\cite{AW}.

\bigskip 

\textbf{Acknowledgement.} We would like to thank Dr. G. Dirr for a
valuable comment.


\end{document}